\documentstyle[11pt,amscd,amssymb]{article}
\newtheorem{thm}{Theorem}[section]
\newtheorem{lem}[thm]{Lemma}

\newtheorem{defi}[thm]{Definition}

\newcommand{\inv}{^{-1}}
\newcommand{\iso}{\stackrel{\sim}{\longrightarrow}}
\newcommand{\tim}{^\times}

\newcommand{\Z}{{\mathbb{Z}}}

\newcommand{\Q}{{\mathbb{Q}}}
\newcommand{\PP}{{\mathbb{P}}}

\newcommand{\cod}{\mbox{codim}}
\newcommand{\cO}{{\cal{O}}}

\newcommand{\ov}{\overline}
\newcommand{\dist}{\mbox{dist}}
\newcommand{\rk}{\mbox{rk}}

\newcommand{\om}{\Omega}
\newcommand{\di}{\mbox{div}}

\begin{document}
\title{Non-Archimedean intersection indices on projective spaces and the Bruhat-Tits building for $PGL$}
{ \author{Annette Werner \\ \small Mathematisches Institut, Universit\"at M\"unster, Einsteinstr. 62, D -  48149 M\"unster\\ \small e-mail: werner@math.uni-muenster.de}}
\date{}
\maketitle
\centerline{\bf Abstract}

Inspired by Manin's approach towards a geometric interpretation of Arakelov theory at infinity, we interpret in this paper non-Archimedean local intersection numbers of linear cycles in $\PP^{n-1}$ with the combinatorial geometry of the Bruhat-Tits building associated to $PGL(n)$. 

\vspace*{1.5cm}

\begin{center}
\LARGE Indices d'intersection non archim\'ediens sur les espaces \\ projectifs et l'immeuble de Bruhat-Tits de $PGL$
\end{center}

\vspace*{0.5cm}

\centerline{\bf R\'esum\'e}

Nous interpr\'etons les nombres d'intersection non archim\'ediens des cycles lin\'eaires sur $\PP^{n-1}$ avec la g\'eom\'etrie combinatoire de l'immeuble de Bruhat-Tits associ\'e \`a $PGL(n)$. Ces r\'esultats sont inspir\'es par le travail 
de Manin donnant une interpr\'etation g\'eo\-m\'etrique de la th\'eorie d'Arakelov  aux places infinies pour les courbes. 

\vspace*{1.5cm}
\centerline{{\bf MSC}(2000): 14C17, 14M15, 20E42}

\centerline{{\bf Keywords: }intersection theory, projective spaces, Bruhat-Tits building}

\newpage
\section{Introduction}

In this paper we interpret non-Archimedean local intersection numbers of linear cycles in $\PP^{n-1}$ with the combinatorial geometry of the Bruhat-Tits building associated to $PGL(n)$. 

The ultimate motivation behind these results is to associate to a projective space a differential geometric object playing the role of a model at infinity in the sense of Manin, who constructed in \cite{ma} such an object for curves. 

A first step in this direction is to look for a geometric interpretation of non-Archi\-me\-dean intersection numbers which has an Archimedean analogue. 
It turns out that the Bruhat-Tits building for $PGL$ is a good candidate for such a geometric framework. 

The goal of the present paper is to express non-Archimedean intersection numbers in terms of the building.
In another work (see \cite{we}) we construct the desired ``model at infinity'' for projective spaces as an Archimedean analogue of the building. Besides, we use results of the present paper to derive parallel geometric formulas for Archimedean and non-Archimedean Arakelov intersection numbers. 

Let us now describe our main results. 
We denote by $X$ the Bruhat-Tits building associated to the group $G = PGL(V)$, where $V$ is an $n$-dimensional vector space over a non-Archimedean local field $K$ of characteristic $0$. The vertices in $X$ correspond to the homothety classes $\{M\}$ of $R$-lattices $M$ in $V$, where $R$ is the ring of integers in $K$.

We fix a lattice $M$ in $V$, which induces a projective space $\PP(M)$ over $R$, and consider $d$ linear cycles on $\PP(M)$ intersecting properly in a cycle of dimension 0. These cycles are equal to projective spaces $\PP(N_i)$ for split $R$-submodules $N_i$ of $M$. We put $L_j = \cap_{i \neq j} N_i$, and define $F$ as the following set of vertices in $X$:
\[F = \{\{\pi^{k_1} L_1 + \ldots + \pi^{k_d} L_d\}: k_1,\ldots, k_d \in \Z\}.\]
In Theorem 5.1 we express Serre's intersection number $<\PP(N_1), \ldots, \PP(N_d)>$ of our linear cycles in $\PP(M)$ as:
\[ <\PP(N_1), \ldots, \PP(N_d)>= \dist(\{M\}, F),\]
where $\dist$ is what we call the combinatorial distance function in $X$, i.e. the minimal length of a path consisting of 1-simplices connecting $\{M\}$ with a vertex in $F$. 

In the case of $\PP^1$, this result specializes to a formula in \cite{ma}.

Moreover we investigate the case of several linear cycles meeting properly in a cycle of higher dimension. In this case, of course, we no longer have an intersection number, but we can nevertheless describe the intersection cycle (see Theorem 5.2). It consists of one component coming from the generic fibre, which appears with multiplicity 1, and one component coming from the special fibre, appearing with multiplicity $\dist(\{M\}, F)$, where $F$ is defined in a similar way as in the previous result. 

{\bf Acknowledgements: } I thank Ch. Deninger, G. Kings, K. K\"unnemann, E. Landvogt, Y.I. Manin, P. Schneider,  E. de Shalit and M. Strauch for useful and inspiring discussions. I am also grateful to the Max-Planck-Institut f\"ur Mathematik in Bonn for financial support and the stimulating atmosphere during the early stages of this project.

\section{Intersection theory}
In this section we will list the definitions and results from  intersection theory which are needed later on, thereby fixing our notation. 

Let $\Omega$ be a scheme, of finite type and separated over a regular ring. By $Z^p(\om)$ we denote the codimension $p$ cycles on $\om$, i.e. the free abelian group on the set of integral (i.e. irreducible and reduced) closed subschemes of codimension $p$ . We write $CH^p(\om)$ for the quotient of $Z^p(\om)$ after the group generated by the principal cycles $\di(f)$ for rational functions $f \neq0$ on a codimension $p-1$ integral closed subscheme, see \cite{fu}, 1.3.

If $T \subset \om $ is a closed subset, let $Z^p_T(\om)$ denote the free abelian group on the set of codimension $p$ integral closed subschemes of $\om$, which are contained in $T$ and $CH^p_T(\om)$ the Chow group of cycles supported on 
$T$, i.e. $Z^p_T(\om)$ modulo the subgroup generated by all $\di(f)$ for rational functions $f \neq 0$ on some codimension  $p-1$ integral closed subscheme of $\om$ which is contained in $T$.  If $T = \emptyset$, we put $CH^p_T(\om) = 0$. 

Now we will briefly recall Serre's intersection pairing. It is defined on any smooth scheme $\Omega$ over a discrete valuation ring by \cite{se}, V-32. Two closed, integral subschemes $X$ and $Y$ of $\Omega$ meet properly if  for every irreducible component $W$ of $X \cap Y$ we have
$\cod(X) + \cod(Y) = \cod(W)$. By \cite{se} we always have the inequality $''\geq''$. If $X$ and $Y$ meet properly, then Serre defines an intersection index $i_W(X,Y)$ of $X$ and $Y$ along $W$ by higher Tor functors, see \cite{se}, V-21.

Let $W$ be an irreducible component of $X \cap Y$, let  $\cO_{\Omega,w}$ be the local ring at the generic point $w$ of $W$, and let $p_X, p_Y$ be the ideals in $\cO_{\Omega, w}$ corresponding to $X$ and $Y$. If $X$ and $Y$ are locally Cohen-Macaulay, then
\[i_W(X,Y) = l_{\cO_{\Omega,w}}(\cO_{\Omega,w} / (p_X+p_Y))\]
by \cite{se}, p. V-20.

We define the intersection  cycle of properly intersecting $X$ and $Y$ by $X\cdot Y =\sum_{W} i_W(X,Y) W$ where the sum runs over all irreducible components of $X \cap Y$.  We can continue this product linearly to arbitrary cycles $X$ and $Y$ meeting properly on $\Omega$, which means that any irreducible component of $X$ meets all the  irreducible components of $Y$ properly. If $X_1,\ldots, X_r$ are $r$ closed integral subschemes meeting properly, i.e. so that every irreducible component $W$ of $\cap_i  X_i$ satisfies $\cod W = \sum \cod X_i$, then the intersection of $X_1, \ldots X_r$ is defined inductively: $X_1 \cdot \ldots \cdot X_r = (\ldots((X_1 \cdot X_2) \cdot X_3)\ldots)\cdot X_r$. Again we can extend this product linearly to arbitrary cycles $X_1, \ldots, X_r$ meeting properly. 

Gillet and Soul\'e have defined an arithmetic intersection pairing for arithmetic Chow groups in \cite{giso}. The local contributions of this pairing at the finite places (in the smooth case) can be described as follows (see \cite{giso}, 4.5.1):

Assume that $\pi: \om \rightarrow S= \mbox{Spec} \; R$ is smooth, separated and of finite type over the discrete valuation ring $R$ and that $\om$ is irreducible. Let $Z_p(\om/S)$ be the free abelian group on the closed integral subschemes $Y \subset  \om$ of relative dimension $p$ over $S$. Here the relative dimension  $\dim_S(Y)$ 
of $Y$ over $S$ is defined as 
\[ \dim_S Y = \mbox{transcendence degree of  }k(Y)\mbox{ over }k(T)- \cod_S(T),\] 
where $T$ is the closure of $\pi(Y)$ in $S$ and $k(Y)$, $k(T)$ are the function fields.
This relative dimension has the property that 
\[ \dim_S(\om) = \dim_S(Y) + \cod_\om(Y)\]
for all closed integral subschemes $Y$ of $\om$ (see \cite{fu}, Lemma 20.1).

Let $CH_p(\om/S)$ be  $Z_p(\om/S)$ modulo rational equivalence. For all closed subschemes $T \subset \om$ we have $CH_T^p(\om) = CH_{d-p}(T/S)$, where $d$ is the relative dimension of $\om$ over $S$. 

For closed subschemes $Y$ and $Z$ we can define a pairing
\[CH_Y^p(\om) \times CH_Z^q(\om) \longrightarrow CH^{p+q}_{Y \cap Z}(\om)\]
as follows: It suffices to define a pairing
\[CH_{d-p}(Y/S) \times CH_{d-q}(Z/S) \longrightarrow CH_{d-p-q}(Y \cap Z/S).\]
Let $V \subset Y$ and $W \subset Z$ be integral closed subschemes. Then we define a cycle in $Y \times Z$ as follows:
\[V \otimes W= \left\{ \begin{array}{ll} 
~ 0, & \mbox{ if } V \mbox{ and } W \mbox{ are contained in the closed fibre},\\
~[V \times_S W], & \mbox{ otherwise}\\
\end{array} \right. \]
By \cite{fu}, Proposition 20.2, this induces a pairing
\[CH_{d-p}(Y/S) \times CH_{d-q}(Z/S) \longrightarrow CH_{2d-p-q}(Y \times_S Z/S).\]

Since $\om$ is smooth over $S$, the diagonal embedding $\Delta: \om \rightarrow \om \times_S \om$ is a regular embedding, and we have Fulton's Gysin map (see \cite{fu}, \S 6 and \S 20):
\[ \Delta^!: CH_{2d-p-q}(Y\times_S Z /S) \longrightarrow CH_{d-p-q}(Y \cap Z/S).\]
Hence we get the desired pairing.

If $V$ and $W$ meet properly,  their image under this pairing coincides with the image of Serre's intersection pairing in $CH^{p+q}_{Y \cap Z}(\om)$ by \cite{fu}, 7.1.2 and 20.2.2.

In particular, if $Y$ and $Z$ are irreducible with $p= \cod_\om (Y)$ and $q = \cod_\om(Z)$, the image of $(Y,Z)$ via $CH_Y^p(\om) \times CH_Z^q(\om) \rightarrow CH^{p+q}_{Y \cap Z}(X)$ yields an intersection class
\[ Y \cdot Z = \Delta^! (Y \otimes Z) \in CH^{p+q}_{Y \cap Z}(\om).\]

If one cycle is given by a Cartier divisor, say $Y = D$, then we have a different description of the image of $D \cdot Z$ in $CH^{p+q}_Z(\om)$ (see \cite{fu}, 
8.1.1, 20.2.1 and 6.1c):
It is equal to the class of $E$, where $E$ is any Weil divisor on $Z$ induced by a Cartier divisor whose line bundle is isomorphic to $j^\ast \cO(D)$. Here $j: Z \hookrightarrow \om$ is the embedding of $Z$ into $\om $ and $\cO(D)$ is the line bundle on $\om$ corresponding to the class of $D$. In particular, if $D$ and $D'$ are linear equivalent divisors on $\om$, the images of $D \cdot Z$ and $D' \cdot Z$ in $CH^{p+q}_Z(\om)$ coincide.

If we have several irreducible cycles $Y_1$, $Y_2,\ldots ,Y_r$ in $\om$ of codimensions $p_1$, $p_2,\ldots, p_r$  we can define inductively an intersection class $Y_1 \cdot \ldots \cdot Y_r \in CH^{p_1 + \ldots + p_r}_{\cap Y_i}(\om)$.

Let $k$ be the residue field of $R$. We denote by $\deg$ the degree map for 0-cycles in the special fibre $\om_k$ of $\om$, i.e. for all $z = \sum n_P P \in Z^d(\om_k)$ we put $\deg z = \sum n_P \; [k(P):k]$, where $k(P)$ is the residue field of $P$.

Assume additionally that $\om$ is proper over $S$, and let  $Y \in Z^p(\om)$ and $Z \in Z^q(\om)$ be two irreducible closed subschemes such that $p+q=d+1$ which intersect properly on the generic fibre of $\om$. This means that their generic fibres are disjoint, so that $Y \cap Z$ is contained in the special fibre $\om_k$ of $\om$. Hence  we can define an intersection number
\[<Y,Z> = \deg( Y \cdot Z),\]
where we take the degree of the image of  $Y \cdot Z \in CH^{d+1}_{Y \cap Z}(\om)$ in $CH^d(\om_k)$.
Similarly, if $Y_i \in Z^{p_i}(\om)$ for $i=1, \ldots,r$ are prime cycles with $\sum p_i = d+1$ which meet properly on the generic fibre, we put $<Y_1, \ldots ,Y_r> = \deg(Y_1 \cdot  \ldots \cdot Y_r)$.

\section{Hyperplanes}
Throughout this paper we denote by $K$ a finite extension of $\Q_p$, by $R$ its valuation ring and by $k$ the residue class field. Besides, $v$ is the valuation map, normalized so that it maps a prime element to $1$. We write $q$ for the cardinality of the residue class field, and we normalize the absolute value on $K$ so that $|x|= q^{-v(x)}$. 

Besides, we fix an $n$-dimensional vector space $V$ over $K$.

Let $\PP(V)=\mbox{Proj Sym}V^\ast$ be the projective space corresponding to $V$, where $V^\ast$ is the linear dual of $V$. Every non-zero linear subspace $W$ of $V$ defines an integral (i.e. irreducible and reduced) closed subscheme $\PP(W)=\mbox{Proj Sym}W^\ast \hookrightarrow  \PP(V)$ of codimension $n - \dim W$. 

By ``$R$-lattice in $V$'' we always mean an $R$-lattice in $V$ of full rank. Every $R$-lattice $M$ in $V$ defines a model $\PP(M) = \mbox{Proj Sym}_R (M^\ast)$ of $\PP(V)$ over $R$, where $M^\ast$ is the $R$-linear dual of $M$. If the lattices $M$ and $N$ differ by multiplication by some $\lambda \in K\tim$ then the corresponding isomorphism $\PP(M) \iso \PP(N)$ induces the identity on the generic fibre.

Throughout this paper we call a submodule $N$ of $M$ split, if the exact sequence $0 \rightarrow N \rightarrow M \rightarrow M/N \rightarrow 0$ is split, i.e. if $M/N$ is free (or, equivalently, torsion free).
Every split $R$-submodule $N$ of $M$ defines a closed subscheme  $\PP(N)=\mbox{Proj Sym}N^\ast \hookrightarrow \PP(M)$. 

\begin{sloppypar}
\begin{lem}
For every split $R$-submodule $N$ of $M$, the closed subscheme $\PP(N)=\mbox{Proj Sym}\, N^\ast$ of $\PP(M)$ is integral, and has codimension $n - \mbox{\rm rk} N$.
\end{lem}
\end{sloppypar}

{\bf Proof: } This follows from dualizing the sequence $ 0 \rightarrow N \rightarrow M \rightarrow M/N \rightarrow 0$.\hfill$\Box$

The cycles in $\PP(M)$ induced by split submodules are called linear, and linear cycles of codimension $1$ are called hyperplanes.  The homogeneous prime ideal corresponding to the linear cycle $\PP(N) \hookrightarrow \PP(M)$ is generated by a base of $(M/N)^\ast$ (regarded in $M^\ast$). In particular, it is generated by homogeneous elements of degree one. 

Now fix a lattice $M$ in $V$ and an $R$-basis $x_1,\ldots, x_n$ of $M$. Let $B$ be a matrix in $GL(n,R)$ which we regard as an endomorphism of $M$ via our fixed basis. Then $B$ induces an automorphism (which we also denote by $B$) of $\PP(M)$. The following lemma can be proven easily: 

\begin{lem}
If the hyperplane $H$ in $\PP(M)$  is given by the linear homogeneous element 
\[f = \sum_{j=1}^n a_j x_j^\ast \in M^\ast,\]
where $x_1^\ast,\ldots,x_n^\ast$ is the dual basis of $x_1,\ldots,x_n$, then $B(H)$ is given by the homogeneous element $\sum_{j=1}^n b_j x_j^\ast$ where 
\[\left( \begin{array}{c} b_1 \\ \vdots \\ b_n \end{array} \right)
= ~^t B\inv \left( \begin{array}{c} a_1\\ \vdots \\ a_n  \end{array} 
\right)  .\]
\end{lem}

Now we need an easy matrix lemma. We call a quadratic matrix a permutation matrix if it contains exactly one entry $1$ in every line and column, and if all other entries are equal to zero.

\begin{lem} Let $A = (a_{ij})$ be an $(n \times n)$-matrix over $R$. Then there exist elements $C$ and $D$ in $GL(n,R)$, where $D$ is a permutation matrix, such that the matrix $CAD = (b_{ij})_{i,j}$ is upper triangular with
\[v(b_{11}) \leq \ldots \leq v(b_{nn}) \mbox{ and}\]
\[v(b_{ii}) \leq v(b_{ij}) \mbox{ for all } i \leq j.\]
\end{lem}
{\bf Proof:} We move a coefficient with minimal valuation in the upper left corner and eliminate the other entries in the first column. This can be repeated until our matrix is upper triangular.\hfill$\Box$

A crucial step for our geometric formulas for intersection indices (to be proven in section 5) is the expression of the intersection number of $n$ hyperplanes in terms of their equations. We can do this for any $n$ hyperplanes $H_1, \ldots, H_n$ in $\PP(M)$ such that their generic fibres $H_{1K}, \ldots, H_{nK}$ meet properly on $\PP(V)$. 

\begin{thm} Let $M$ be a lattice in $V$. We fix a basis $x_1,\ldots,x_n$ of $M$, and denote by $x_1^\ast, \ldots, x_n^\ast \in M^\ast$ the dual basis. Let $H_1,\ldots,H_n$ be hyperplanes in $\PP(M)$ which intersect properly on the generic fibre. 
Let $f_i=\sum_j a_{ji} x_j^\ast \in M^\ast$ be a linear  homogeneous element generating the ideal  corresponding to $H_i$ and put $A=(a_{kl})_{k,l}$. Then we have the following formula for the intersection number of $H_1, \ldots, H_n$:
\[<H_1,\ldots, H_n>= v(\det A).\]
\end{thm}

{\bf Proof: }By Lemma 3.3 we find  some $C \in \mbox{GL}(n,R)$ and a permutation matrix $D$ such that $CAD=B = (b_{ij}) $ is upper triangular and satisfies the inequalities
\[v(b_{11}) \leq \ldots \leq v(b_{nn}) \mbox{ and}\]
\[v(b_{ii}) \leq v(b_{ij}) \mbox{ for all } i \leq j.\]
There is a permutation $\sigma$ of $\{1,\ldots n\}$ such that $AD$ is the coefficient matrix for the hyperplanes $H_{\sigma(1)}, \ldots, H_{\sigma(n)}$. 
By Lemma 3.2, the linear element $\sum_{j \leq i} b_{ji} x_j^\ast$ corresponds to the hyperplane $~^t C\inv (H_{\sigma(i)})$ for all $i = 1,\ldots, n$. Now $<H_1,\ldots , H_n> = $ \linebreak[3] $<~^t C\inv (H_{\sigma(1)}) , \ldots ,~^tC\inv (H_{\sigma(n)})>$ and $v(\det A) = v( \det B)$. Hence we can assume that $A$ is upper triangular with $v(a_{11}) \leq \ldots  \leq v(a_{nn})$ and $v(a_{ii}) \leq v(a_{ij})$ for $i \leq j$. 

We can assume that $\cap H_i \neq \emptyset$, since otherwise our claim is trivial. 

For all $a \in R$ we denote by $\overline{a}$ its image in $k$. The reduction $(H_i)_k$ of $H_i$ corresponds to the homogeneous ideal generated by $\overline{f_i}=\sum\overline{a_{ji}} x^\ast_j$ in $\PP(M_k^\ast)$. Let $\overline{A}$ be the matrix $(\overline{a_{ij}})$. Now $\cap (H_{i})_k$ is the linear cycle corresponding to the subspace $L_k \subset M_k$ which is equal to $\mbox{ker}~^t\overline{A}$ via the identification of $M_k$ with $k^n$ given by the reductions of $x_1, \ldots, x_n$.

We will first assume that $H_1,\ldots,H_n$ meet properly on the whole of $\PP(M)$. 
Hence their intersection consists of one point in the special fibre, and $v(a_{11}) = \ldots = v(a_{n-1  n-1})=0$. We can therefore assume that 
$a_{11} = \ldots = a_{n-1 n-1} = 1$.
Hence for all $k \leq n-1$ the homogeneous ideal in $\mbox{\rm Sym } M^\ast$ generated by $f_1 = x_1^\ast, f_2 = a_{12} x_1^\ast + x_2^\ast, \ldots, f_k = \sum_{j=1}^{k-1} a_{jk} x_j^\ast + x_k^\ast$ is equal to the homogeneous ideal generated by $x_1^\ast,\ldots,x_k^\ast$. Now we can compute Serre's intersection index as follows:

Note that all $H_i$ are isomorphic to $\PP_R^{n-2}$, hence they are locally Cohen-Macaulay (even regular). Assume that $n >2$. The closed subset $H_1 \cap H_2$ of $\PP(M)$ is given by the homogeneous ideal $(f_1,f_2) = (x_1^\ast, x_2^\ast)$; the corresponding reduced closed subscheme is the linear cycle $\PP(N) \hookrightarrow \PP(M)$ for $N = (M^\ast / R x_1^\ast + R x_2^\ast)^\ast$. Hence $W = H_1 \cap H_2$ is a prime cycle. We have $i_W(H_1,H_2) =1$, hence $H_1\cdot H_2 = W$. Besides, $W$ is a projective space over $R$, hence also locally Cohen-Macaulay.

The same argument (if $n >3$) implies that $H_1 \cdot H_2 \cdot H_3$ is equal to the cycle given by the irreducible subset $H_1 \cap H_2 \cap H_3$. Finally we find that $H_1 \cdot H_2 \cdot \ldots \cdot H_{n-1} =W$ where $W$ is the prime cycle corresponding to the homogeneous ideal $(x_1^\ast, \ldots, x_{n-1}^\ast)$. Now it is easy to calculate $<H_1,\ldots, H_n> = v(a_{nn}) = v( \det A)$, which proves our claim in the case of proper intersection.

Hence we can now assume that there is an $s < n-1$ such that $v(a_{11}) = \ldots = v(a_{ss}) = 0$ and $v(a_{kk}) >0$ if $k>s$. We write $l_i = v(a_{ii})$, and we can again assume that $a_{11} = \ldots = a_{ss} = 1$.  For all $m= 1,\ldots, n-1$ let $Y_m$ be the cycle corresponding to the integral subscheme given by the homogeneous ideal $(x_1^\ast, \ldots, x_m^\ast)$ of $\mbox{Sym } M^\ast$, and let $Z_m$ be the cycle corresponding to the integral subscheme given by $(x_1^\ast,\ldots,x_{m-1}^\ast,\pi)$, where $\pi$ is a fixed prime element in $R$.

Assume that $m$ is a number with $s < m < n$. Then the intersection of $Y_{m-1}$ and $H_m$ has two irreducible components, namely $Y_m$ and $Z_m$. Since both have codimension $m$, the cycles $Y_{m-1}$ and $H_m$ meet properly. We want to calculate $i_{Y_m}(Y_{m-1}, H_m)$ and $i_{Z_m}(Y_{m-1}, H_m)$. Let $y$ respectively $z$ be the generic points of $Y_m$ respectively $Z_m$. Since $m <n$, they are both contained in $U = \{x_n^\ast \neq 0\}$. We write $A = \cO_{\PP(M), y}$. Then $A = R[y_1,\ldots,y_{n-1}]_{(y_1,\ldots,y_m)}$ with $y_i = x_i^\ast / x_n^\ast$. Since $Y_{m-1}$ and $H_m$ are locally Cohen-Macaulay, we have
\[ i_{Y_m}(Y_{m-1}, H_m) = l_{A}(A/(y_1,\ldots,y_{m-1},\pi^{l_m} y_m)).\]
As $\pi$ is a unit in $A$, this is equal to $l_A(A/(y_1,\ldots,y_m)) = 1$.

Similarly, we put $B = \cO_{\PP(M), z}$, hence $B = R[y_1,\ldots,y_{n-1}]_{(y_1,\ldots,y_{m-1},\pi)}$, and  we get
\[i_{Z_m}(Y_{m-1}, H_m) = l_{B}(B/(y_1,\ldots,y_{m-1},\pi^{l_m} y_m)).\]
Here $y_{m}$ is a unit in $B$, hence this length is equal to 
\[l_B(B/(y_1,\ldots,y_{m-1},\pi^{l_m})) = l_{(R[y_m,\ldots, y_{n-1}]/(\pi^{l_m}))_{(\pi)}}((R[y_m,\ldots, y_{n-1}]/(\pi^{l_m}))_{(\pi)}) = l_m,\]
since the only ideals in $(R[y_m,\ldots, y_{n-1}]/(\pi^{l_m}))_{(\pi)}$ are $0$ and $(\pi^k)$ for $0 < k < l_m$ and these are all distinct.

Now we will prove by induction that $H_1 \cdot \ldots \cdot H_m$ is equal to the class of  $Y_m + (\sum_{i \leq m} l_i ) Z_m$ in $CH^m_{\cap_{i \leq m }H_i}(\PP(M))$ for all $m \leq n-1$. Since $H_1$ is irreducible, we have $l_1=0$ and $H_1 =Y_1$, which is our claim for $m =1$. 

Now we come to the induction step. Assume that our claim holds for some $m$ with $1 \leq m \leq n-2$. If $m \leq s$, we can move to $m+1$ using the above calculations of intersection indices. Let us now assume that $m > s$, hence that $l_m$ and $l_{m+1}$ are strictly positive.
Then $H_{m+1}$ meets $Y_m$ properly in the components $Y_{m+1}$ and $Z_{m+1}$ and we can calculate $Y_m \cdot H_{m+1}$ via Serre's intersection: $Y_m \cdot H_{m+1}$ is induced by the cycle  $Y_{m+1 } + l_{m+1} Z_{m+1}$ by our previous calculations. 

Note that $f_{m+1}$ is contained in $(x_1^\ast, \ldots, x_{m-1}^\ast, \pi)$, since both $l_m$ and $l_{m+1}$ are strictly positive. Hence $Z_m$ is contained in $H_{m+1}$. 
Now we can determine the intersection 
\[Z_m \cdot H_{m+1} \in CH_{Z_m \cap H_{m+1}}^{m+1}(\PP(M)) = CH_{Z_m}^{m+1}(\PP(M))\]
by the recipe we described in section 2 for intersections with divisors. Since $H_{m+1}$ is linearly equivalent to the hyperplane $H'$ given by the ideal $(x_{m}^\ast)$, we find that $Z_m \cdot H_{m+1}$ is equal to the image of $Z_m \cdot H'$ in $CH_{Z_m}^{m+1}(\PP(M))$. Now $Z_m \cap H' = Z_{m+1}$ and $Z_m $ and $H'$ meet properly in this irreducible set with $i_{Z_{m+1}}(Z_m,H')=1$, which implies that $Z_m \cdot H'$ is induced by the cycle  $Z_{m+1}$.

Altogether we find that $H_1 \cdot \ldots \cdot H_{m+1}$ is the image of $Y_{m+1} + (l_1 + \ldots + l_{m+1}) Z_{m+1}$ in $CH^{m+1}_{\cap_{i \leq m+1} H_i}(\PP(M))$, which finishes the proof of our claim.

We know now that $H_1 \cdot \ldots \cdot H_{n-1}$ is the image of $Y_{n-1} + (l_1 + \ldots + l_{n-1}) Z_{n-1}$ in $CH^{n-1}_{\cap_{i \leq n-1} H_i}(\PP(M))$. Since $s < n-1$, we have $ l_{n-1}>0$ and $l_n >0$. Now $Y_{n-1}$ meets $H_n$ properly, hence we can apply our result in the case of hyperplanes meeting properly on the whole of $\PP(M)$ and find 
$\deg (Y_{n-1} \cdot H_n) = l_n$.

Besides, we have $Z_{n-1} \cap H_n = Z_{n-1}$, and as in our induction step we can show that $Z_{n-1} \cdot H_n$ is the image of $Z_n$ in $CH^{n}_{Z_{n-1}}(\PP(M))$. Since $Z_n$ is a $k$-rational point in the special fibre, we get $\deg ( Z_{n+1} \cdot H_n) = 1$.

Altogether we find that
\[ <H_1,\ldots,H_n> = l_1 + \ldots+ l_n = v(\det A),\]
whence our claim.\hfill$\Box$

\section{The Bruhat-Tits building for $PGL$}
We denote by $X$ the Bruhat-Tits building corresponding to the group $G = PGL(V)$ (see \cite{brti1}). 
$X$ is a metric space with a continuous $G$-action and a simplicial structure. For our purposes, we can think of it as the geometric realization of the following simplicial complex:  We call two lattices in $V$ equivalent, if they differ by a factor in $K\tim$, and we write $\{M\}$ for the equivalence class of the lattice $M$. Two different lattice classes $\{M'\}$ and $\{N'\}$ are called adjacent, if there are representatives $M$ and $N$ of $\{M'\}$ and $\{N'\}$ such that 
\[\pi N \subset M \subset N.\]
This relation defines a flag complex, namely  the simplicial complex whose vertex set is the set of all classes $\{M\}$, and whose simplices are the sets of pairwise adjacent lattice classes. Note that it carries a natural $G$-action.

If $n=2$, then $X$ is an infinite regular tree, with $q+1$ edges meeting in every vertex. 

The building $X$ is the union of its apartments, which correspond to the maximal split tori in $G$. We can describe them as follows: For every decomposition $V = \bigoplus_{1 \leq i \leq n} L_i$  of $V$ in one-dimensional subspaces $L_i$ generated by some vector $v_i$ we define an apartment as the subcomplex of $X$ given by all lattices $M$ which can be diagonalized with respect to our decomposition, i.e. $M = \sum_{i = 1}^n  R \pi^{k_i} v_i$ for some integers $k_i$. 

\begin{defi}The combinatorial distance $\dist(x,y)$ of two vertices $x$ and $y$ in $X$ is defined as 
\begin{eqnarray*}
\dist(x,y) &= &\min \{ k: \mbox{ there are vertices }x= x_0, x_1 ,\ldots, x_k = y,\\
~  &~&\quad \quad \mbox{so that }x_i \mbox{ and }x_{i+1}\mbox{ are adjacent for all }i=0,\ldots, k-1.\}.
\end{eqnarray*}
\end{defi}
Hence $\dist$ is the minimal number of $1$-simplices forming a path between $x$ and $y$. Note that $\dist$ is in general not proportional to the metric on $X$. 

\begin{lem}
Let $x = \{M\}$ and $y = \{L\}$ be two vertices in $X$, and define
\begin{eqnarray*}
s& = &\min\{k: \pi^k L \subset M\} \quad \mbox{and}\\
r & = & \max \{k : M \subset \pi^k L \}.
\end{eqnarray*}
Then we have $\dist(x,y) = s-r$.
\end{lem}
{\bf Proof: }Note that the term on the right hand side is independent of the choice of a representative of the lattice classes. 

Put $d = \dist(x,y)$. Then we find lattices $M = M_0, M_1, \ldots ,M_d$ such that $M_d = \alpha L$ for some $\alpha \in K\tim$ and such that
\[ \pi^d M_d \subset \pi^{d-1} M_{d-1} \subset\ldots \subset \pi M_1 \subset M_0 \subset M_1 \subset \ldots \subset M_d.\]
Hence $\pi^d \alpha L \subset M$, which implies $s \leq d + v(\alpha)$, and $M \subset \alpha L$, which implies $r \geq v(\alpha)$. Altogether we find that $s-r \leq d$.

Let us now show that also $s-r \geq d$ is true. We have by definition $\pi^s L \subset M$ and $M \subset \pi^r L$. Put $L' = \pi^r L$. By the invariant factor theorem, we find an $R$-basis $w_1,\ldots ,w_n$ of $L'$, such that $M = \sum \pi^{k_i} R w_i$ for some integers $k_i$. Since $\pi^{s-r} L' \subset M \subset L'$, all $k_i$ are between $0$ and $s-r$. 

Now put $l_{ij} = \max\{0, k_j-i\}$ for all $i \in \{0,\ldots, s-r\}$ and all $j \in \{1,\ldots, n\}$. Then we define for $i= 0, \ldots, s-r$
\[ M_i = \pi^{l_{i1}} R w_1+ \ldots + \pi^{l_{in}}R w_n.\]
Note that $M_0 = M$, and $M_{s-r} = L'$. We have for all $i= 0, \ldots, s-r-1$ the inclusions $\pi M_{i+1} \subset M_i \subset M_{i+1}$.
Since either $\{M_i \}= \{M_{i+1}\}$, or $\{M_i\}$ and $ \{M_{i+1}\}$ are adjacent, we found a chain of adjacent lattices of length $\leq s-r$ connecting $\{M\}$ and $\{L\}$, which implies our claim. \hfill$\Box$

\section{Intersection indices via combinatorial geometry}
Let us fix  a lattice $M$ in $V$. We will now interpret non-Archimedean intersection  numbers of linear cycles on $\PP(M)$ with the combinatorial geometry of $X$.

Take some $d \geq2$ and let  $\PP(N_1), \ldots, \PP(N_d)$ be $d$ linear cycles on $\PP(M)$ with 
\[\sum_{i=1}^d \cod \, \PP(N_i) = n,\] 
which meet properly on $\PP(M)$. Then $N_1,\ldots, N_d$ are split submodules of $M$ of rank $r_1,\ldots, r_d$ satisfying $\sum_{i=1}^d r_i = (d-1) n$ by 3.1. We will always assume that $r_i \neq n$. For all $j = 1, \ldots, d$ we put $L_j = \cap_{i \neq j} N_i \subset M$. Since all $M/N_i$ are torsion free, the same holds for  $M/L_j$, so that $L_j$ is a split submodule of $M$.

Now let $F$ be set of all vertices in $X$ of the form 
\[\{\pi^{k_1} L_1 \oplus \ldots \oplus \pi^{k_d} L_d\} \quad \mbox{for some} \quad k_1,\ldots, k_d \in \Z.\]
(Alternatively, we can also work with the convex hull of $F$ in $X$.) Note that the intersection $\cap_{i=1}^d N_{i}$ is zero, as $\PP(N_1), \ldots, \PP(N_d)$ do not meet on the generic fibre. Therefore the sum $\pi^{k_1} L_1 \oplus \ldots \oplus \pi^{k_d} L_d$ is direct. On the other hand, an easy calculation shows that $\dim L_{jK} \geq n-r_j$. Since $\sum_j \dim L_{jK} \geq \sum_j (n-r_j)=n$, we must have $\dim L_{jK} = n - r_j$ and $V= \bigoplus _j L_{jK}$, so that $\pi^{k_1} L_1 \oplus \ldots \oplus \pi^{k_d} L_d$ is indeed a lattice of full rank in $V$. 

Obviously, $F$ is the set of vertices contained in either a full apartment or an intersection of affine hyperplanes in some apartment. (Hence its convex hull is either an apartment of an intersection of affine hyperplanes.)

Let us describe $F$ in two special cases: 

{\bf i)} The case of $n$ hyperplanes: 

First note that $n$ hyperplanes $H_{1 K}, \ldots, H_{n K}$ in  the generic fibre $\PP(V)$ intersecting properly, i.e. not at all, in $\PP(V)$, define an apartment $A(H_{1 K}, \ldots, H_{n K})$ in $X$ as follows: For all $i$ the hyperplane $H_{iK}$ is the linear cycle corresponding to an $(n-1)$-dimensional subspace $W_i$ of $V$. Since the $H_{iK}$ intersect properly, $\cap W_i$ is equal to $0$. For all $j=1,\ldots,n$ put $U_j=\cap_{i \neq j} W_i$. Then all $U_j$ are one-dimensional and $V = \oplus U_j$. We denote the apartment corresponding to this decomposition (see section 4) by $A(H_{1 K}, \ldots, H_{n K})$.

If $H_1,\ldots,H_n$ are hyperplanes in $\PP(M)$ which intersect properly on the whole of $\PP(M)$, then the subset $F$ is just the apartment $A(H_{1 K}, \ldots H_{n K})$ corresponding to the generic fibres.

{\bf ii)} The case of two cycles: 

Let $N_1$ and $N_2$ be two  non-trivial split submodules of $M$ of rank $p$ respectively $q= n-p$, such that the corresponding linear cycles $Z_1 = \PP(N_1)$ and $Z_2 = \PP(N_2)$ meet properly on $\PP(M)$. Then the corresponding subset $F$ is
\[ \{\{N_1 + \pi^k N_2\}: k \in \Z\}.\]
In fact, one can show that $F$ defines a doubly infinite geodesic in the building $X$ whose boundary points (with respect to the Borel-Serre compactification of $X$) are the parabolics induced  by the vector spaces $W_i = N_i \otimes_R K$ for $i= 1$ and $2$. 

We will now show that the intersection number of $d$ linear cycles meeting properly on $\PP(M)$ is the combinatorial distance of the set $F$ to the lattice $M$.

If $n$ is equal to $2$ (i.e. $X$ is the Bruhat-Tits tree associated to $PGL(2,K)$), and we consider the case of two hyperplanes in $\PP(M)$, then the ensuing formula is due to Manin (see \cite{ma}, p. 232). Note that in this setting our special cases i) and ii) from above coincide, and the subset $F$ we are dealing with is a geodesic in the Bruhat-Tits tree.

\begin{thm}
The intersection number of $\PP(N_1), \ldots, \PP(N_d)$ on $\PP(M)$ can be expressed as follows with the combinatorial geometry of $X$:

\[<\PP(N_1), \ldots, \PP(N_d)> = \dist(\{M\},F),\]

where $\dist$ denotes the combinatorial distance in $X$, and $\{M\}$ is the vertex in $X$ defined by our fixed lattice $M$.
\end{thm}
{\bf Proof: } Recall that $r_i$ is the rank of $N_i$ and that the numbers $m_i = n-r_i$ satisfy $\sum_{i=1}^d m_i = n$. Besides, put $n_0 = 0 $ and $n_i = m_1 + \ldots + m_i$ for $i = 1,\ldots,d$. Then $n = n_d$.

For all $i = 1,\ldots, d$ let $g_{n_{i-1}+1},\ldots, g_{n_i}$ be a $R$-basis of $(M/N_i)^\ast \hookrightarrow M^\ast$. The elements $g_{n_{i-1}+1}, \ldots, g_{n_i}$ generate the homogeneous prime ideal corresponding to $\PP(N_i) \hookrightarrow \PP(M)$. Besides, fix a basis $y_1,\ldots, y_n$ of $M$, and let $A'=(a'_{ij})$ be the coordinate matrix of $g_1,\ldots, g_n$ with respect to the dual basis $y_1^\ast,\ldots, y_n^\ast$, i.e. 
\[ g_j = \sum_{i=1}^n a'_{ij} y_i^\ast.\]

Now choose an element $a'_{ij}$ with $j \leq n_1$ such that $v(a'_{ij})$ is minimal among the entries of the first $n_1$ columns. We remove the $j$-th column and insert it before the first one. Then we switch rows so that this element sits in the upper left corner, and we perform some elementary row operations to eliminate $a'_{21}, \ldots, a'_{n1}$. Among the $j$ between $2 $ and $n_1$ and among the $i \geq 2$ we choose an entry of minimal valuation. Again, after removing the corresponding column and putting it between the first and the second one, and after elementary row operations involving only the last $n-1$ rows we can assume that $a'_{i2} = 0$ for all $i >2$. We continue in this way until we reach the $n_1$-th column. Then $a'_{ij}$ is zero for $j \leq n_1$ and $i>j$, and the upper left corner satisfies the following divisibility conditions:
\begin{eqnarray*}
v(a'_{ii})\leq v(a'_{i+1 i+1})  \quad \mbox{for } i \leq n_1-1 \quad \mbox{and }\\
v(a'_{ii}) \leq v(a'_{ij}) \quad \mbox{for }  i < j \leq n_1 . \quad\quad\quad~ 
\end{eqnarray*}
Now choose an  element $a'_{ij}$ with $i,j \geq n_{1}+1$ and $j \leq n_2$ such that $v(a'_{ij})$ is minimal among the entries $a'_{ij}$ with $i \geq n_1 +1$ and $n_1 +1 \leq j \leq n_2$. As before, we permute the columns corresponding to $g_{n_{1}+1}, \ldots, g_{n_2}$ and perform elementary row operations involving only the rows with index bigger or equal to $n_1+1$ to achieve $a'_{ij} = 0$ for $j \leq n_2$ and $i>j$. 

We continue this process until we worked our way through the whole matrix. We see that after permuting the equations $g_{n_{i-1}+1 }, \ldots, g_{n_i}$ corresponding to each of the  $N_1, \ldots, N_d$,  and after switching to another basis of $M$ (in order to take care of the row operations) we can assume that our coordinate matrix $A'$ is upper triangular. 

Let us denote the reduction of elements in $R$ or $M$ or of $R$-matrices by overlining. Besides we use the following notation for submatrices: For any $n\times n$-matrix $D= (d_{kl})$ and any $ i \leq j \leq n$ we write 
$D(i\rightarrow j)$ for the $n \times (j-i+1)$-submatrix consisting of the columns $i, i+1,\ldots, j$. Similarly we write $D(i \downarrow j)$ for the $(j-i+1) \times n$-submatrix consisting of the rows $i,i+1,\ldots,j$. By $D(i_1 \rightarrow j_1, i_2 \rightarrow j_2)$ we mean the submatrix where we take columns $i_1, \ldots, j_1$ followed by columns $j_1,\ldots, j_2$. Besides, $D(i_1 \downarrow j_1)(i_2 \rightarrow j_2)$ is the submatrix consisting of all entries $d_{kl}$ with $i_1 \leq k \leq j_1$ and $i_2 \leq l \leq j_2$. 

For $n\times n$-matrices we will furthermore abbreviate $D_{ij} = D(n_{i-1}+1 \downarrow n_i)(n_{j-1}+1\rightarrow n_j)$, if $i$ and $j$ are $\leq d$. If we divide $D$ into rectangular submatrices according to our partition $n = m_1 + \ldots + m_d$, then $D_{ij}$ is the rectangle in position $(ij)$. 

Since $g_{n_{i-1}+1}, \ldots, g_{n_i}$ form a basis of  the split $R$-submodule $(M/N_i)^\ast$ of $M^\ast$, their reductions generate a vector space of rank $m_i$ over $k$, so that the coordinate matrix 
$\overline{A}'(n_{i-1}+1 \rightarrow n_i)$ has full rank $m_i$. 

Since $\PP(N_1), \ldots, \PP(N_d)$ meet properly, their intersection is empty or zero-dimen\-sional and contained in the special fibre. Hence $\cap_i \PP({N}_{ik})$ is empty or has dimension zero, so that $\cap_{i} {N}_{ik}$ has dimension $\leq 1$. Since ${N_{ik}}$ is equal to the kernel of $(\ov{g}_{n_{i-1}+1}, \ldots, \ov{g}_{n_i})$, this means that $rk \ov{A}' \geq n-1$. 

The intersection of $\PP(N_1), \ldots, \PP(N_d)$ is empty iff $\rk \, \ov{A}' = n$. In this case all elements on the diagonal of $A'$ are units. After switching to another basis of $M$, we can therefore assume that $A'$ is diagonal. Then $M = L_1 + \ldots + L_d$, so that $\{M\}$ is actually contained in $F$, and our formula holds. 

Hence we only have to deal with the case $\rk \, \ov{A}'= n-1$.
Here the intersection of $\PP(N_1),\ldots, \PP(N_d)$ is not empty. As $v(\det A') > 0$, we find  indices $p \leq q \leq d$ with
\begin{eqnarray*}
v(\det A'_{11}) = \ldots = v(\det A'_{p-1 p-1}) = 0,\\
v(\det A'_{pp}) > 0, \quad v(\det A'_{qq}) >0\quad \mbox{and}\\
v(\det A'_{q+1 q+1} ) = \ldots = v(\det A'_{dd}) = 0.
\end{eqnarray*}

Note that $p \geq 2$, since $\ov{A}'(1 \rightarrow n_1)$, and thus $\ov{A}'_{11}$ has full rank. 
All elements on the diagonal of $A'_{11},\ldots, A'_{p-1 p-1}$, $A'_{q+1 q+1},\ldots, A'_{dd}$ are units. 
After performing some elementary row operations we can assume that $A'_{ij} = 0$, if $i < j \leq p-1$ or if $j \geq q+1$ and $i<j$. Note that the divisibility conditions in $A'_{qq}$ imply that all elements on the diagonal except possibly the last one are units. 
Therefore we can eliminate all entries $a'_{ij}$ in $A'$ such that $1 \leq i \leq n_{q-1}$ and 
$n_{q-1}+1 \leq j \leq n_q-1$ by elementary row operations. Hence after switching to another basis $y_1,\ldots, y_n$ of $M$ we may assume that our coordinate matrix $A'$ contains zeroes above $A_{22}', \ldots, A_{p-1 p-1}'$ and $A'_{q+1,q+1}, \ldots A'_{dd}$, and that all columns above $A'_{qq}$ are zero except possibly the last one.

Now recall that for $j = 1,\ldots, d$ the module $L_j =\cap_{i \neq j} N_i$ is a split submodule of $M$ of rank $m_j$. For each $j=1,\ldots,  d$ choose an $R$-basis
\[w_{n_{j-1}+1}, \ldots , w_{n_j}\]
of $L_j$. We denote by $B'$ the transpose of the coordinate matrix of $w_1, \ldots, w_n$ with respect to $y_1,\ldots, y_n$, i.e. $B'$ is the matrix $B' = (b'_{ij})_{ij}$ so that
\[ w_i = \sum_{j=1} ^n b'_{ij} y_j.\]
For all $j = 1,\ldots, d$ there are matrices 
$C_j, D_j \in GL(m_j, R)$ such that 
$
C_j B'_{jj} D_j $ is a diagonal matrix with entries $\beta_{ii}$ such that $v(\beta_{11}) \leq \ldots \leq v(\beta_{m_j m_j})$. Let $C$ (respectively $D$) be the matrix with diagonal components $C_1, \ldots, C_d$ (respectively $D_1, \ldots, D_d$), 
and define $B = C B' D$. 

Then $D$ is the transpose of the transition matrix from a base $x_1, \ldots, x_n$ of $M$ to our base $y_1,\ldots,y_n$. If we put
\[v_i = \sum_{j=1}^n b_{ij} x_j,\]
then for all $h \leq d$ the elements $v_{n_{h-1}+1}, \ldots, v_{n_h}$ form a basis of $L_h$.
The matrix $B$ is by definition the transpose of the transition matrix from $x_1,\ldots, x_n$ to $v_1,\ldots, v_n$. 

Now $D\inv$ is the transition matrix from the dual basis $x_1^\ast, \ldots, x_n^\ast$ of $M^\ast$ to $y_1^\ast, \ldots, y_n^\ast$. Hence the coordinate matrix of $g_1,\ldots, g_n$ with respect to $x_1^\ast, \ldots, x_n^\ast$ is equal to $A'' = D\inv A'$.
The matrix $A''$ has the property that if a block $A'_{ij}$ is zero or zero up to the last column, then the same holds for $A''_{ij}$. Besides, $v(\det A''_{jj}) = v(\det A'_{jj})$. 
However, the diagonal blocks
$A''_{11}, \ldots, A''_{d d}$ may not be upper triangular any longer.

After permuting the first $n_1 = m_1 $ columns we may assume that $v(a''_{n_1 n_1}) \leq v(a''_{n_1 j})$ for all $j = 1,\ldots, n_1$. By a series of elementary column operations we can eliminate $a''_{n_1 1},\ldots, a''_{n_1 n_1-1}$. Now we permute the first $n_1 -1$ columns to achieve $v(a''_{n_1-1 j}) \leq v(a''_{n_1-1  n_1-1})$ for all $j = 1,\ldots, n_1-1$, and we clear out 
$a''_{n_1-1 1},\ldots, a''_{n_1-1 n_1-2}$. Note that these column operations affect only the first $n_1-1$ columns, hence the first $n_1-1$ elements in the $n_1$-th row remain zero. We repeat this process until $A''_{11}$ is upper triangular. These column operations amount to passing to another basis $f_1,\ldots, f_{n_1}$ of $(M/N_1)^\ast$

Now we work on $A''(n_1+1 \rightarrow n_2)$. First we switch columns inside this block to achieve $v(a''_{n_2 j}) \leq v(a''_{n_2  n_2})$ for all $j = n_1 + 1,\ldots, n_2-1$, and we clear out 
$a''_{n_2 n_1 +1},\ldots, a''_{n_2 n_2 -1}$, then we eliminate $a''_{n_2-1 n_1 +1},\ldots, a''_{n_2-1 n_2 -2}$, and so on, until $A''_{22}$ is upper triangular. We do this with block after block  until $A''$ is upper triangular. 

Since the columns were transformed block by block, we can find for each $i = 1,\ldots, d$ a new $R$-basis
\[f_{n_{i-1}+1}, \ldots, f_{n_i}\]
of $(M/N_i)^\ast$, whose coordinate matrix $A$ with respect to the basis $x_1^\ast,\ldots, x_n^\ast$ is upper triangular and has the property that the columns above $A_{22}, \ldots, A_{p-1 p-1}$ and above $A_{q+1 q+1} \ldots, A_{dd}$ are zero. Besides, we still have $v(\det A_{pp}) >0$ and $v(\det A_{qq}) >0$. Since there exists a matrix $D \in GL(m_q, R)$ such that 
\[A(n_{q-1}+1 \rightarrow n_q)(1 \downarrow n_{q-1}) = A''(n_{q-1}+1 \rightarrow n_q)(1 \downarrow n_{q-1}) D,\]
we still know that $\rk \,  \ov{A}(n_{q-1}+1 \rightarrow n_q)(1 \downarrow n_{q-1}) \leq 1$. 

Let $H_j$ be the hyperplane given by the linear homogeneous element $f_j$. Then $H_{n_{i-1}+1}, \ldots, H_{n_i}$ intersect properly and $\PP(N_i) = H_{n_{i-1}+1}\cdot  \ldots \cdot H_{n_i}$. Hence by 3.4 we have 
\[<\PP(N_1),\ldots, \PP(N_d)> = v(\det A).\]

Besides note that
\[f_j(v_i)= \sum_k a_{kj} b_{ik},\]
so that the entries of $BA$ are equal to  $(f_j(v_i))_{ij}$.

Now by definition, $N_i$ lies in the kernel of all $f_{n_{i-1}+1}, \ldots, f_{n_i}$. Since $L_j$ is contained in $N_i$ for all $i \neq j$, this implies that the blocks $(BA)_{ij}$ are equal to $0$  for $i \neq j$. Since all $B_{jj}$ are diagonal matrices,  $B$ is therefore upper triangular. Besides we have
$B_{ij} = 0 $, if $i <j \leq p-1$, and if $i<j$ and $j >q$. 

Since $v_{n_{j-1}+1}, \ldots, v_{n_j}$ are a basis of the split submodule $L_j \subset M$, their reductions are still linear independent, hence
the rank of $\ov{B}(n_{j-1}+1 \downarrow n_j)$ is equal to $m_j$. We will now show that only the last diagonal element $b_{n_j n_j}$ in the block $B_{jj}$ may not be a unit. Indeed, assume that $v(b_{ii})>0$ for some $n_{j-1}+1 \leq i \leq n_j-1$. By the divisibility conditions along the diagonal of $B_{jj}$ we have $v(b_{n_j-1 n_j-1})>0$ and $v(b_{n_j n_j}) >0$. For all $l = n_{j-1}+1, \ldots, n_j$ we find $f_l(v_{n_j-1}) = b_{n_j-1 n_j-1} a_{n_j-1 l }$, so that $\ov{f_l(v_{n_j-1})} = 0$. Similarly, $\ov{f_l(v_{n_j})}=0$, so that $\ov{v_{n_j-1}}$ and $\ov{v_{n_j}}$ are contained in $N_{jk}$.  Since $L_j$ is contained in all $N_i$ for $i \neq j$, the elements $\ov{v_{n_j-1}}$ and $\ov{v_{n_j}}$ lie also in $N_{ik}$ for all $i \neq j$. Hence we found two linear independent vectors in $\cap_i N_{ik}$, which contradicts our assumption that $\PP(N_1), \ldots, \PP(N_d)$ meet properly. 

We will now show the following claim

(1) There exists an index $t \leq p-1$ such that $\ov{A}(1 \rightarrow n_{t-1}, n_t+1 \rightarrow n)$ has full rank.

Let us first assume that $p<q$.
Note that $\rk \, \ov{A}(n_{q-1}+1 \downarrow n_q)$ is strictly smaller than $m_q$, since $v(\det(A_{qq})) >0$. Since $\rk\,  \ov{A} = n-1$, the matrix $\ov{A}(n_{p-1}+1 \downarrow n_{q-1})$ must have full rank $(n_{q-1} - n_{p-1})$. Disregarding the zeros, we find that also the $(n_{q-1} -n_{p-1}) \times (n_q - n_{p-1})$-matrix 
\[\ov{A}(n_{p-1}+1 \downarrow n_{q-1})(n_{p-1} +1 \rightarrow n_q)\] 
must have full rank $(n_{q-1} - n_{p-1})$. This matrix consists of the two vertical chunks
\[\ov{A}(n_{p-1}+1 \downarrow n_{q-1})(n_{p-1} +1 \rightarrow n_{q-1}) \mbox{ and }\ov{A}(n_{p-1}+1 \downarrow n_{q-1})(n_{q-1} +1 \rightarrow n_q).\]
We know that the second chunk $\ov{A}(n_{p-1}+1 \downarrow n_{q-1})(n_{q-1} +1 \rightarrow n_q)$ has rank $\leq 1$. Hence the rank of the first chunk
\[\Omega:=\ov{A}(n_{p-1}+1 \downarrow n_{q-1})(n_{p-1} +1 \rightarrow n_{q-1}) \]
must be $\geq (n_{q-1} - n_{p-1} -1)$. 

Besides, as $v(\det A_{pp}) >0$, there must be an element on the diagonal of $A_{pp}$ which has positive valuation, i.e. there exists an index $n_{p-1} + 1 \leq i \leq n_p$ such that $\ov{a}_{ii}=0$. Since $\Omega$ has rank $\geq (n_{q-1} - n_{p-1} -1)$, the upper left corner $\ov{A}(n_{p-1}+1 \downarrow i-1)(n_{p-1} +1 \rightarrow i)$ must have full rank $( i -1 -n_{p-1} )$. Now recall that $\ov{A}(n_{p-1}+1 \rightarrow n_p)$ has full rank, hence the first $i -n_{p-1}$ columns of this matrix, namely $\ov{A}(n_{p-1}+1 \rightarrow i)$, also have full rank $i - n_{p-1}$. Put
\[\lambda_j = \ov{A}(n_{p-1}+1 \rightarrow i)(j\downarrow j),\]
which is just the $j$-th row  of this matrix. We have seen that $\lambda_{n_{p-1}+1}, \ldots, \lambda_{i-1}$ are linear independent, hence there exists a row $\lambda_{j_0}$ with $j_0 \leq n_{p-1}$ so that $\lambda_{j_0}, \lambda_{n_{p-1}+1}, \ldots,$\linebreak[3]$\lambda_{i-1}$ is a full linear independent set of rows in $\ov{A}(n_{p-1}+1 \rightarrow i)$.

Now let $t \leq p-1$ be the index of the block in which $\lambda_{j_0}$ lies, i.e. $n_{t-1} + 1 \leq j_0 \leq n_t$. 
We want to show that $\ov{A}(1 \rightarrow n_{t-1}, n_t +1 \rightarrow n)$ has full rank. It obviously suffices to show that $
\ov{A}(n_{t}+1 \rightarrow n)(n_{t-1}+1 \downarrow n)$ has full rank. 
Removing the first block of rows $\ov{A}(n_{t}+1 \rightarrow n)(n_{t-1}+1 \downarrow n_{t})$ and putting it before the block  $\ov{A}(n_{t}+1  \rightarrow n)(n_{p-1}+1 \downarrow n_p)$, we see that it remains to show that the matrix
\begin{eqnarray*}
\Lambda = \left(
{\footnotesize
\begin{array}{ccc}
\ov{A}_{tp} & \ldots& \ov{A}_{td}   \\
\ov{A}_{pp} & \ldots& \ov{A}_{pd} \\
& \ddots  & \vdots  \\
0& & \ov{A}_{dd}  \\
\end{array}
}
\right).
\end{eqnarray*}

has full rank $n-n_{p-1}$. Recall that we fixed an element $a_{ii}$ on the diagonal of $A_{pp}$ with positive valuation. Now
\begin{eqnarray*}
\Lambda =  \left(
{\small 
\begin{array}{ccrlccc}
\ov{a}_{n_{t-1}+1, n_{p-1}+1} & \ldots  &\ov{a}_{n_{t-1}+1,i} |& & &  \\
&  \vdots &  | && \ast& \\
\ov{a}_{n_t, n_{p-1}+1} & \ldots& \ov{a}_{n_t i}|&&&\\
\ov{a}_{n_{p-1}+1, n_{p-1}+1} &\ldots & \ov{a}_{n_{p-1}+1, i}| &&&\\
& \vdots &| &&&\\
\ov{a}_{i-1, n_{p-1}+1}& \ldots  &\ov{a}_{i-1 i}| && &    \\ 
--- & -- & ---| &---& -- & --\\
&   & | &\ov{a}_{i i+1} & \ldots &\ov{a}_{i n} \\
& & | & \ov{a}_{i+1 i+1} &\ldots  &  \ov{a}_{i+1 n}\\
& 0&| & & \ddots &   \vdots\\
& &| & & &\ov{a}_{nn} \\ 
\end{array}
}
\right).
\end{eqnarray*}
Hence in the upper left corner of $\Lambda$ we find the rows $\lambda_j$ for $j=n_{t-1}+1, \ldots, n_t$ and $j = n_{p-1}+1, \ldots, i-1$, and we know that there are $(i-n_{p-1})$ linear independent ones among them. So the upper left corner of $\Lambda$ has full rank $i-n_{p-1}$. The lower right corner (which is a $(n-i+1) \times (n-i)$-matrix) has the same rank as the $(n-i+1) \times n$-matrix $\ov{A}(i \downarrow n)$. Since $\rk (\ov{A}) = n-1$, the rank of the lower right corner must therefore be bigger or equal to $n-i$ which means that this corner has also full rank. Hence $\Lambda$ has full rank. This finishes the proof of claim (1) in the case $p <q$. 

If $p$ is equal to $q$, similar arguments can be used to prove $(1)$.

Now we claim that we can calculate our intersection number as follows: For the index $t$ we found in (1) we have 

(2)\hfill$<\PP(N_1),\ldots, \PP(N_d)> = v(b_{n_t n_t})$.\hfill~

Put $w_1= v_{n_{t-1}+1}, \ldots, w_{m_t} = v_{n_t}$, which is a basis of $L_t = \cap_{i \neq t} N_i$. Since $L_t$ is a split submodule of $M$, we can complete this basis to a basis $w_1, \ldots, w_n$ of $M$. 
Since $f_j = \sum_{i=1}^n f_j(w_i) w_i^\ast$, the coordinate matrix of $f_1,\ldots, f_n$ with respect to the dual basis $w_1^\ast,\ldots, w_n^\ast$ is equal to $(f_j(w_i))_{i,j}$. Using 3.4, we can therefore calculate the intersection number as follows:
\[ <\PP(N_1),\ldots, \PP(N_d)> = v(\det(f_j (w_i)_{i,j=1,\ldots,n})).\]
Now $f_j(w_i)=0$ for $i=1,\ldots, m_t$ and $j \notin J_t:=\{n_{t-1}+1, \ldots, n_t\}$. Hence after permuting columns the matrix $(f_j (w_i))_{i,j=1,\ldots,n}$ looks like this $\left( \begin{array}{cc} \ast & 0 \\ \ast & \ast \end{array} \right)$, and we get:
\[v(\det(f_j(w_i)_{i,j=1,\ldots,n})) = v(\det(f_j (w_i)_{i=1,\ldots, m_t, j \in J_t })) + v(\det(f_j (w_i)_{i= m_t+1,\ldots, n, j \notin J_t})).\]
Note that $(f_j(w_i))_{i=1,\ldots, m_t, j \in J_t } = B_{tt} A_{tt}$. As $t $ is strictly smaller than $p$, the determinant of $A_{tt}$ is a unit. Besides, as we have shown above, only the last element on the diagonal of $B_{tt}$ may not be a unit,  so that we can calculate the first term as follows:
\[v(\det(f_j(w_i)_{i=1,\ldots, m_t, j \in J_t })) = v(b_{n_t  n_t}).\]
It remains to show that the second term is zero, hence that $(\ov{f_j(w_i)})_{i= m_t+1,\ldots, n, j \notin J_t}$ is an invertible $(n-m_t ) \times (n-m_t) $-matrix over $k$. Since $f_j(w_i) = 0$ for $i=1, \ldots, m_t$ and $j \notin J_t$ we can as well show that the matrix $(\ov{f_j(w_i)})_{i= 1,\ldots, n, j =1, j \notin J_t}$  has full rank $n-m_t$. In order to prove this we may change the base of $M$ and show that $\ov{A}(1 \rightarrow n_{t-1}, n_t+1 \rightarrow n) = (\ov{f_j(x_i)})_{i= 1,\ldots, n, j =1, j \notin J_t}$ has full rank $n-m_t$, which was done in (1).

Let us fix some $h= 1, \ldots, d$. We will first show by induction that for all $i  \geq n_h$ we have

(3)\hfill $b_{n_h n_h} \mbox{ divides } b_{n_h i } \prod_{j= n_h+1}^i a_{jj},$\hfill~

where the empty product is equal to $1$. For $i = n_h$ this is trivial. So let us suppose our claim is true for all $i$ with $n_h \leq i < i_0$ for some $i_0 \geq n_h+1$. Then 
\[0 = f_{i_0}(v_{n_h}) = \sum_{i= n_h}^{i_0}  b_{n_h i} a_{i i_0}.\]
We multiply by $\prod_{j= n_h+1}^{i_0-1}a_{jj}$ and get
\[ 0 = \sum_{i= n_h}^{i_0-1} b_{n_h i}(\prod_{j= n_h+1}^{i_0-1}a_{jj})a_{i i_0} + b_{n_h i_0 }(\prod_{j= n_h+1}^{i_0-1}a_{jj})a_{i_0 i_0}, \]
hence the induction hypothesis implies our claim for $i_0$.

Note that (3) implies that $v(b_{n_t n_t}) = \max_{i=1,\ldots, n}v( b_{ii})$. Indeed, if $i$ is not equal to $n_h$ for some $h= 1, \ldots, d$, then $b_{ii}$ is a unit. Since $\ov{B}$ contains no zero lines, for every $n_h$ there must be an index $i \geq n_h$ such that $v(b_{n_h i}) = 0$. By (3) we find that 
\[v(b_{n_h n_h}) \leq \sum_{j=n_h+1} ^i v( a_{jj}) \leq v(\det A).\]
But $v(\det A)$ is equal to the intersection number $<\PP(N_1),\ldots, \PP(N_d)>$, hence equal to $ v(b_{n_t n_t})$ by $(2)$. 

Now we can prove 

(4) \hfill $\dist(\{M\}, F) = v(b_{n_t n_t}).$\hfill~

We take a lattice class $\{L\}$ corresponding to a vertex in $F$, i.e. $L = \pi^{k_1} L_1 + \ldots + \pi^{k_{d-1}} L_{d-1} + L_d$ for some integers $k_1,\ldots, k_{d-1}$. Then $\min\{k: \pi^k L \subset M\}\geq -k_t$ and $\max \{k : M \subset \pi^k L \} \leq -k_t -v(b_{n_t n_t})$, so that 4.2 implies $\dist(\{M\}, \{L\}) \geq v(b_{n_t n_t})$.
Hence we are done if we show that 
\[\dist(\{M\},\{ L_1+ \ldots+ L_d\}) \leq v(b_{n_t n_t}).\]

Recall that $BA$ is the matrix consisting of the diagonal blocks $B_{11} A_{11}, \ldots, B_{dd} A_{dd}$. Now for all $h = 1, \ldots, d$ the determinant of $A_{hh}$ is $\prod_{j = n_{h-1}+1}^{n_h} a_{jj}$, and $v(\det B_{hh}) = v(b_{n_h n_h})$. Applying again (3) and the fact that in every line in $B$ there must be a unit, we find 
\[v(b_{n_h n_h}) \leq v\left(\prod_{j=n_h+1}^n a_{jj}\right),\]
so that  
\[v(\det B_{hh} A_{hh}) \leq v\left(\prod_{j = n_h+1}^n a_{jj}\right)+ v\left(\prod_{j=n_{h-1}+1}^{n_h} a_{jj}\right)  \leq v(\det A) = v(b_{n_t n_t}),\]
hence $\pi^{ v(b_{n_t n_t})} (A_{hh}\inv B_{hh}\inv)$ has $R$-coefficients for all $h$. 
Therefore $\pi^{v(b_{n_t n_t}) }B\inv$ has $R$-coefficients, which implies $ M \subset \pi^{-v(b_{n_t n_t})} (L_1+ \ldots + L_d)$, as $~^t\! B M = L_1 + \ldots + L_d$. Now it is easy to see that indeed $\dist(\{M\},\{ L_1 + \ldots + L_d\}) \leq v(b_{n_t n_t})$.\hfill $\Box$

One may wonder if our assumption that the linear cycles meet properly on the whole of $\PP(M)$ is really necessary. 
In fact, in 3.4 we have  proven a formula for the intersection number of $n$ hyperplanes meeting properly {\em only }on $\PP(V)$. Hence it is tempting to try and use this as a starting point to generalize Theorem 5.1 at least to the case of $n$ hyperplanes $H_1, \ldots, H_n$ in $\PP(M)$ meeting properly only on the generic fibre $\PP(V)$. Unfortunately, the geometric expression $\dist(\{M\}, A(H_{1K},\ldots, H_{nK})) $ does in general no longer coincide with the intersection number $<H_1, \ldots, H_n>$ in this case. In fact, one can show as in the proof of 5.1 that for $H_1,\ldots, H_n$ with coordinate matrix $A = (a_{ij})$ (upper triangular and subject to $v(a_{ii}) \leq v(a_{i+1 i+1})$ and $v(a_{ii}) \leq v(a_{ij})$ for $i <j$)
we have $\dist(\{M\}, A(H_{1K},\ldots, H_{nK})) = v(a_{nn})$, whereas the intersection number \mbox{$<H_1, \ldots, H_n>$} $ = \sum_{i=1}^n v(a_{ii})$ (see 3.4) may be bigger.

We will conclude this paper by generalizing Theorem 5.1 in another direction. 
Namely, let us consider the case of $d \geq 2$ linear cycles on $\PP(M)$ meeting properly in a cycle of dimension bigger than zero. So let $N_1, \ldots, N_d$ be split submodules of $M$ of rank $r_1\ldots, r_d$, so that $\PP(N_1), \ldots, \PP(N_d)$ are linear cycles meeting properly on $\PP(M)$, i.e. any irreducible component in the intersection $\cap_i \PP(N_i)$ has codimension 
\[ \sum_i \cod \,\PP(N_i) = dn -\sum_{i=1}^d r_i.\]
Note that only the case $\sum_i \cod \,\PP(N_i) \leq n$ is interesting, since otherwise $\cap \PP(N_i)$ must be empty. 
Since we already dealt with the case $\sum_i \cod \, \PP(N_i) = n$ in 5.1, let us now assume that $\sum_i \cod \, \PP(N_i) < n$. 
Put $r_0 = n - (dn -\sum_i r_i) = \sum_i r_i - (d-1)n >0$, and let $L_0$ be the intersection $L_0 = \cap_{i=1}^d N_i$. Since all $M/N_i$ are torsion free, the same holds for $M/L_0$, so that $L_0$ is a split submodule of $M$. Hence $\PP(L_0)$ is a linear cycle contained in the intersection $\cap_i \PP(N_i)$, so it has codimension $\geq n-r_0$, which implies that $\rk \, L_0 \leq r_0$. On the other hand, we can calculate
\[ \rk \, L_0 = \dim (\cap N_{iK}) \geq \sum_{i=1}^d \dim (N_{iK}) - (d-1) n = \sum_{i=1}^d r_i - (d-1)n = r_0,\]
so that we find $\rk \, L_0 = r_0$. 

For $j = 1, \ldots d$ we put $L_j = \cap_{i \neq j} N_i \subset M$. Then $L_0 $ is also a split submodule of $L_j$, hence there exists a free $R$-module $L_j'$ with $L_j = L_0 \oplus  L_j'$. Let $F$ be the set of vertices in $X$ of the form
\[ \{\pi^{k_1} L_1' \oplus \ldots \oplus \pi^{k_d} L_d' \oplus \pi^{k_{d+1}} L_0 \}\quad \mbox{for some }\quad k_1,\ldots, k_{d+1} \in \Z. \]
Of course, $F$ depends on the choice of $L_1', \ldots, L_d'$.
It is easy to see  that the sum $
\pi^{k_1} L_1' +\ldots + \pi^{k_d} L_d' + \pi^{k_{d+1}} L_0$ is indeed direct. Since 
\[\dim L_{iK} \geq \sum_{j \neq i} \dim N_{jK} - (d-2)n = \sum_{j \neq i} r_j - (d-2)n = n-r_i+r_0\]
we find that 
\begin{eqnarray*}
\lefteqn{ \sum_{i=1}^d \dim L'_{iK} + \dim L_{0K} = }\\
& & \sum_{i=1}^d \dim L_{iK} - (d-1) r_0 \geq dn -\sum_{i=1}^d r_i +  r_0 = n,
\end{eqnarray*}
so that $\pi^{k_1} L_1' \oplus \ldots \oplus \pi^{k_d} L_d' \oplus \pi^{k_{d+1}} L_0$ is indeed a lattice of full rank in $V$. Besides we see that  $\rk \, L'_i  = n-r_i$. 

Then the strategy of the proof of 5.1 can be modified to prove the following result:

\begin{thm}
Let $\PP(N_1), \ldots, \PP(N_d)$ be linear cycles meeting properly on $\PP(M)$ such that $\sum_{i=1}^d \cod \; \PP(N_i) <n$.
For any choice of $F$ as above, we can describe the Serre intersection cycle $\PP(N_1) \cdot \ldots \cdot \PP(N_d)$ as follows:
\[\PP(N_1) \cdot \ldots \cdot \PP(N_d) =\ov{\PP(\cap_{i=1}^d N_{iK})} + \dist(\{M\}, F) \; \PP(\cap_{i=1}^d N_{ik}).\]
Here the first cycle $\ov{\PP(\cap_{i=1}^d N_{iK})}$ is the closure in $\PP(M)$ of the linear cycle $\PP(\cap_{i=1}^d N_{iK})$  on $\PP(V)$, which is just the Serre intersection cycle $\PP(N_1)_K \cdot \ldots \cdot \PP(N_d)_K$ of the generic fibres. The second cycle $\PP(\cap_{i=1}^d N_{ik})$ is the linear cycle on the special fibre corresponding to the subspace $\cap_{i=1}^d N_{ik}$ of $M_k$, which we regard as a cycle on $\PP(M)$. 
\end{thm}

Both our results 5.1 and 5.2 are proven by a complicated series of direct computations once the correct formulas are found. It would be desirable to give a more conceptual proof providing deeper insights in the nature of our formulas and allowing generalizations.

\end{document}